\documentclass[12pt]{article}

\usepackage{latexsym}
\usepackage{amsmath}
\usepackage{amsfonts}
\usepackage{amssymb}

\newcommand{\bean}{\begin{eqnarray}}
\newcommand{\eean}{\end{eqnarray}}
\newcommand{\bea}{\begin{eqnarray*}}
\newcommand{\eea}{\end{eqnarray*}}
\newcommand{\bsa}{\begin{subarray}{c}}
\newcommand{\esa}{\end{subarray}}
\newcommand{\bi}{\begin{itemize}}
\newcommand{\ei}{\end{itemize}}

\newtheorem{lemma}{Lemma}[section]

\newtheorem{propn}[lemma]{Proposition}

\title{\bf Vector-valued modular forms with an unnatural boundary}
\author{Marvin Knopp \\
Department of Mathematics \\
Temple University \\
Philadelphia, Pa. \\
and\\
Geoffrey Mason\thanks{%
Supported by the NSF and NSA} \\
Department of Mathematics, \\
University of California at Santa Cruz, \\
CA 95064.\\
\\
Dedicated to the memory of Leon Ehrenpreis.}
\date{}

\begin{document}
\maketitle

\begin{abstract}
\noindent We characterize all logarithmic, holomorphic vector-valued modular forms
which can be analytically continued to a region strictly larger than the upper half-plane. (MSC2010: 11F12, 11F99.)
\end{abstract}

\section{Introduction}
Set $\Gamma = SL_2(\mathbb{Z})$ be the modular group with standard generators
\begin{eqnarray*}
S = \left(\begin{array}{cc}0 & -1 \\ 1 & 0\end{array}\right), \ T = \left(\begin{array}{cc}1 & 1 \\ 0 & 1\end{array}\right).
\end{eqnarray*} 
Let $\rho: \Gamma \rightarrow GL(p, \mathbb{C})$ be a $p$-dimensional representation of $\Gamma$.
A holomorphic vector-valued modular form of weight $k\in  \mathbb{Z}$ associated to $\rho$ is a holomorphic function $F: \frak{H} \rightarrow \mathbb{C}^p$ defined on the upper half-plane $\frak{H}$ which satisfies
\begin{eqnarray}\label{fe}
F|_k \gamma (\tau) = \rho(\gamma)F(\tau) \ \ \ \ (\gamma \in \Gamma)
\end{eqnarray}
and a growth condition at $\infty$ (see below). As usual, the stroke operator here is defined as
\begin{eqnarray*}
F|_k \gamma (\tau) = (c\tau+d)^{-k}F(\gamma \tau)\ \ (\gamma = \left(\begin{array}{cc} a&b  \\ c&d\end{array}\right)\in \Gamma).
\end{eqnarray*}

We also refer to the pair $(\rho, F)$ as a vector-valued modular form, and call $p$ the \emph{dimension} of $(\rho, F)$. We usually consider $F$ as a vector-valued function\footnote{Superscript t means \emph{transpose} of vectors and matrices}
 $F(\tau) = (f_1(\tau), \hdots, f_p(\tau))^t$ and call the $f_i(\tau)$ the \emph{component functions}
 of $F$. The purpose of the present paper is to investigate whether $F(\tau)$ has a \emph{natural boundary}.

\medskip
If $f(\tau)$ is a nonconstant (scalar) modular form of weight $k$ on a subgroup of finite index in  $\Gamma$, then it is well known that the real axis is a \emph{natural boundary} for $f(\tau)$. That is, there is \emph{no} real number $r$  such that $f(\tau)$ can be analytically continued to a region containing $\frak{H} \cup \{r\}$. In \cite{KM2} the authors extended this result to \emph{normal} vector-valued modular forms on 
$\Gamma$, showing that each nonconstant component of such a modular form
has the real line as natural boundary. Here we study the same question for the larger class of \emph{polynomial} (or 
\emph{logarithmic}) vector-valued modular forms introduced  in \cite{KM3}, where one assumes only that the eigenvalues of $\rho(T)$ have absolute value $1$. This case is more subtle because the existence of a natural boundary is \emph{no longer  true} in general. Here, we say that a vector-valued modular form has the real line as a natural boundary provided that at least one component does.

\medskip
Consider the column vector of polynomials 
\begin{eqnarray*}
C(\tau) = (\tau^{p-1}, \tau^{p-2}, \hdots, 1)^t.
\end{eqnarray*}
It is not hard to see (cf. Section \ref{S2}) that $C(\tau)$ is a vector-valued modular form of weight $1-p$
associated with a representation $\sigma$ equivalent to the
$p-1$st symmetric power $S^{p-1}(\nu)$ of the natural defining representation $\nu$ of $\Gamma$.
Thus for any $p$, $(\sigma, C)$ is a $p$-dimensional holomorphic vector-valued modular form
 which is obviously analytic throughout the complex plane. The main result of the present paper is that this is essentially the only example of a polynomial vector-valued modular form whose natural boundary is not the real line.

 \medskip
 In order to state our results precisely, we need one more definition. Suppose that $(\rho, F), (\rho', F')$ are two vector-valued modular forms of weight $k$ and 
 dimension $p$. We say that
 they are \emph{equivalent} if there is an invertible $p \times p$ matrix $A$ such that
 \begin{eqnarray*}
(\rho, F) = (A\rho' A^{-1}, AF').
\end{eqnarray*}

We give two 
formulation of the main result. As we shall explain, they are essentially equivalent.
 
 \medskip
 \noindent
 {\bf Theorem $1$.} Suppose that the eigenvalues of $\rho(T)$ have absolute value $1$, and let $(\rho, F)$ be a nonzero holomorphic vector-valued modular form of weight $k$ and dimension $p$.  Then the following are equivalent:
 \begin{eqnarray*}
&&(a)\ F(\tau) \ \mbox{does \emph{not} have the real line as a natural boundary},\\
&&(b)\ \mbox{The component functions of $F(\tau)$ span the space of polynomials}\\
&&\ \ \ \ \ \mbox{of degree $l-1$ for some $l \leq p$. Moreover, $k=-l$.}
\end{eqnarray*}

\medskip
\noindent
{\bf Theorem $2$.} Suppose that the eigenvalues of $\rho(T)$ have absolute value $1$, and let $(\rho, F)$ be a holomorphic vector-valued modular form of weight $k$ and dimension $p$. 
Suppose further that the component functions of $F(\tau)$ are \emph{linearly independent}.  Then the following are equivalent:
 \begin{eqnarray*}
&&(a)\ F(\tau) \ \mbox{does \emph{not} have the real line as a natural boundary},\\
&&(b)\ (\rho, F)\ \mbox{is equivalent to}\ (\sigma, C)\ \mbox{and} \ k = 1-p.
\end{eqnarray*}

\medskip
The reader familiar with Eichler cohomology will recognize the space of polynomials in $\tau$ 
(considered as $\Gamma$-module) as a crucial ingredient in that theory. This points to the fact that  Eichler cohomology has close connections to the theory of vector-valued modular forms, connections  that in fact go well beyond the question of natural boundaries that we treat here. The authors hope to return to this subject in the future.

\bigskip
The paper is organized as follows. In Section \ref{S2} we review some background material regarding vector-valued modular forms, discuss the basic example $(\sigma, C)$ introduced above, and explain why Theorems 1 and 2 are equivalent. In Section \ref{S3} we give the proof of the Theorems.

\section{Background}\label{S2}
The space of homogeneous polynomials in variables $X, Y$ is a right $\Gamma$-module 
such that $\gamma= \left(\begin{array}{cc} a&b \\ c&d\end{array}\right) \in \Gamma$ is an algebra automorphism with
\begin{eqnarray*}
\gamma: X \mapsto aX+bY, \ Y \mapsto cX+dY.
\end{eqnarray*}
The subspace of homogeneous polynomials of degree $p-1$ is an irreducible
$\Gamma$-submodule which we denote by $Q_{p-1}$.
The representation of $\Gamma$ that it furnishes is the
$p-1$th symmetric power $S^{p-1}(\nu)$ of the defining representation $\nu$.

\medskip
For $\tau \in \frak{H}$, let $P_{p-1}(\tau)$ be the space of polynomials in $\tau$ of degree
at most $p-1$. Since
\begin{eqnarray*}
\tau^j|_{1-p} \gamma = (c\tau+d)^{p-1} \left( \frac{a\tau+b}{c\tau+d}  \right)^j = (a\tau+b)^j(c\tau+d)^{p-1-j},
\end{eqnarray*}
it follows that $P_{p-1}(\tau)$ is a right $\Gamma$-module with respect to the stroke operator
$|_{1-p}$. Indeed, $P_{p-1}(\tau)$ is isomorphic to $Q_{p-1}$, an isomorphism being given by
\begin{eqnarray*}
X^jY^{p-1-j} \mapsto \tau^j \ \ (0 \leq j \leq p-1).
\end{eqnarray*}

Since $1, \tau, \hdots, \tau^{p-1}$ are linearly independent and span a right $\Gamma$-module with respect to the stroke operator $|_{1-p}$, we know (cf. [KM1], Section 2) that there is a unique representation
$\sigma: \Gamma \rightarrow GL_{p}(\mathbb{C})$ such that
\begin{eqnarray*}
\sigma(\gamma)C(\tau) = C|_{1-p}\gamma (\tau) \ \ (\gamma \in \Gamma).
\end{eqnarray*}
This shows that $(C, \sigma)$ is a vector-valued modular form of weight $1-p$ and that
the representation $\sigma$ is equivalent to $S^{p-1}(\nu)$.

\medskip
We can now explain why Theorems 1 and 2 are equivalent. Assume first that Theorem 1 holds. 
Let $(\rho, F)$ be a holomorphic vector-valued modular form of 
weight $k$ with linearly independent component functions and such that the real line is not a natural boundary for $F(\tau)$. By Theorem 1 the components of $F(\tau)$ span a space of polynomials of degree no greater than $p-1$, and by linear independence  they must span the space 
$P_{p-1}(\tau)$. Moreover, we have $k=1-p$. Now there is an invertible $p \times p$ matrix $A$ such that $AF(\tau)=C(\tau)$,
whence $(\rho, F(\tau))$ is equivalent to $(A\rho A^{-1}, C(\tau))$. As explained above, we necessarily have
$A\rho A^{-1} = \sigma$ in this situation, so that $(\rho, F)$ is equivalent to $(\sigma, C)$. This shows that 
(a) $\Rightarrow$ (b) in Theorem 2, in which case Theorem 2 is true. 

\medskip
Now suppose that Theorem 2 holds, and let $(\rho, F)$ be a nonzero holomorphic vector-valued modular form of dimension $p$ and weight $k$ such that the real line is not a 
natural boundary for $F(\tau)$.
Let $(g_1, \hdots, g_{l})$ be a \emph{basis} for the span of the components of $F$. 
Setting $G = (g_1, \hdots, g_{l})^t$, we again use ([KM1], Section 2) to find a representation
$\alpha: \Gamma \rightarrow GL_{l}(\mathbb{C})$ such that $(\alpha, G)$ is a holomorphic vector-valued modular form of weight $k$. Because the components of $G$ are linearly independent,
 Theorem 2 tells us that they span the space $P_{l-1}(\tau)$ of polynomials of degree at most $l$, and 
 that $k = -l$. Thus the conclusions of Theorem 1(b) hold, and Theorem 1 is true.

\bigskip
To complete this Section we recall (following [KM3])
some facts about the polynomial $q$-expansions which arise as component functions of holomorphic vector-valued modular forms which will be needed for the proof of the main Theorems.

\medskip
Let $(\rho, F)$ be a holomorphic vector-valued modular form of weight $k$.
 Replacing $(\rho, F)$ by an equivalent vector-valued modular form if necessary,
we may,  and shall, 
assume that $\rho(T)$ is in (modified) Jordan canonical form.\footnote{A minor variant of the usual Jordan canonical form. See [KM3] for details.} In passing from $(\rho, F)$  to an equivalent vector-valued modular form, the component functions of $F$ are replaced by linear combinations of  the components, but this will not matter to us. Let the $i$th Jordan block of $\rho(T)$ have size
$m_i$, and label
the corresponding component functions of $F(\tau)$ as
$\varphi^{(i)}_1(\tau), \hdots, \varphi^{(i)}_{m_i}(\tau)$. By [KM3] they have 
polynomial $q$-expansions
\begin{eqnarray}\label{logqexp1}
\varphi_l^{(i)}(\tau) = \sum_{s=0}^{l-1} {\tau \choose s} h^{(i)}_{l-1-s}(\tau) \ \ (1 \leq l \leq m_i), 
\end{eqnarray}
each $h_s^{(i)}(\tau)$ having a left-finite $q$-series
\begin{eqnarray}\label{qexp1}
h_s^{(i)}(\tau) = e^{2\pi i \mu_i \tau}\sum_{n=\nu_i}^{\infty} a_n(s, i)e^{2\pi i n\tau}\ \ (0 \leq s \leq m_i-1,
\ \nu_i \in \mathbb{Z}).
\end{eqnarray}
Here, $e^{2\pi i \mu_i}$ is the eigenvalue of $\rho(T)$ determined by the $i$th block and
$0 \leq \mu_i < 1$. (It is here that we are using the assumption that the eigenvalues of $\rho(T)$ 
have absolute value $1$.) $F(\tau)$ is called \emph{holomorphic} if, for each Jordan block,
each $q$-series $h_s^{(i)}(\tau)$ has only \emph{nonnegative} powers of $q$, i.e.,
$a_n(s, i)=0$ whenever $n+\mu_i<0$.

\section{Proof of the main Theorems}\label{S3}
In this Section we will prove Theorem 2. As we have explained, this is equivalent to Theorem 1.

\medskip
In order to prove Theorem 2, we may replace $(\rho, F)$ by any equivalent vector-valued modular form.
Thus we may, and from now on shall, assume without loss that $\rho(T)$ is in (modified) Jordan canonical form. We assume that $\rho(T)$ has $t$ Jordan blocks, which we may, and shall, further assume are ordered in decreasing size
$M= m_1\geq m_2 \geq \hdots \geq m_t$. Thus $m_1+\hdots+ m_t = p$, and we may speak, with an obvious meaning, of the component functions in a block. The $i$th. block corresponds to an eigenvalue $e^{2\pi i \mu_i}$ of  $\rho(T)$,
and we let the component functions of $F(\tau)$ in that block be as in (\ref{logqexp1}), (\ref{qexp1}).

\bigskip
Let
\begin{eqnarray}\label{abcdmatrix}
\gamma = \left(\begin{array}{cc}a&b \\c&d\end{array}\right) \in \Gamma.
\end{eqnarray}
Because $(\rho, F)$ is a vector-valued modular form of weight $k$ we have
\begin{eqnarray*}
\rho(\gamma)F(\tau) = (c\tau+d)^{-k}F(\gamma \tau).
\end{eqnarray*}
So if $\varphi_u(\tau) = \varphi^{(i)}_v(\tau)$ is the $u$th component of the $i$th block of $F(\tau)$
(so that $u = m_1+\hdots+ m_{i-1}+v$), then
\begin{eqnarray*}
 (c\tau+d)^{-k}\varphi_u(\gamma \tau) = \sum_{j=1}^t \sum_{l=1}^{m_j} \alpha^{(j)}_l \varphi^{(j)}_l(\tau)
\end{eqnarray*}
where $(\hdots, \underbrace{\alpha^{(j)}_1, \hdots, \alpha^{(j)}_{m_j}}_{jth\ block}, \hdots)$ is the $u$th row of $\rho(\gamma)$. Using (\ref{logqexp1}), we obtain
\begin{eqnarray*}
 (c\tau+d)^{-k}\varphi_u(\gamma \tau) &=& \sum_{j=1}^t\sum_{l=1}^{m_j}  \sum_{s=0}^{l-1}
 \alpha^{(j)}_l  {\tau \choose s} h^{(j)}_{l-1-s}(\tau) \notag \\
 &=& \sum_{s=0}^{M-1} {\tau \choose s} \left(\sum_{j=1}^t \sum_{l=s+1}^{m_j} \alpha_l^{(j)}\ h_{l-1-s}^{(j)}(\tau) \right) \notag \\
 &=& \sum_{s=0}^{M-1} {\tau \choose s} \sum_{j=1}^t \left(\alpha_{s+1}^{(j)}h_0^{(j)}(\tau)
 + \sum_{l=s+2}^{m_j} \alpha_l^{(j)}\ h_{l-1-s}^{(j)}(\tau) \right).
  \end{eqnarray*}
(Here, $\alpha_{s+1}^{(j)}=0$ if $s \geq m_j$.)

\medskip
Because the component functions of $(\rho, F)$ are linearly independent, $\varphi_u(\tau)$ is nonzero and the previous display is not identically zero. So
there is a largest integer $B$ in the range $0 \leq B \leq M-1$ such that the summand corresponding to
 ${\tau \choose B}$ does \emph{not} vanish. Now note that
 $\varphi_1^{(j)}(\tau) = h_0^{(j)}(\tau)$. Because the component functions are linearly independent,
  then in particular the  $h_0^{(j)}(\tau)$ are linearly independent, and we can conclude that
 \begin{eqnarray}\label{conc1}
&&\alpha_{s+1}^{(j)} = 0 \ \ (1\leq j \leq t, \ s>B+1), \\
&&\alpha_{B+1}^{(j)}\ \ \mathrm{are\ not\ all\ zero}\ \ (1 \leq j \leq t). \notag
\end{eqnarray}
It follows that
\begin{eqnarray}\label{goodqexp}
&&\ \ \ \ \ \ \ \ \ \ \ \ \ \ \ \ \ \   (c\tau+d)^{-k}\varphi_u(\gamma \tau)  \notag \\
&=& {\tau \choose B}\sum_{j=1}^t 
 \alpha_{B+1}^{(j)}h_0^{(j)}(\tau)+ \sum_{s=0}^{B-1} {\tau \choose s} \sum_{j=1}^t \sum_{l=s+2}^{m_j} \alpha_l^{(j)}\ h_{l-1-s}^{(j)}(\tau)
  \end{eqnarray}
and the first term on the right hand side of (\ref{goodqexp}) is nonzero. Incorporating (\ref{qexp1}), we obtain
\begin{eqnarray}\label{goodqexp1}
&&\ \ \ \ \ \ \ \ \ \ \ \ \    (c\tau+d)^{-k}\varphi_u(\gamma \tau)  \notag \\
&=& {\tau \choose B}\sum_{j=1}^t 
 \alpha_{B+1}^{(j)}e^{2\pi i \mu_j \tau}\sum_{n=\nu_j}^{\infty} a_n(0, j)e^{2\pi i n\tau} \notag \\
 &&+ \sum_{s=0}^{B-1} {\tau \choose s} \sum_{j=1}^t \sum_{l=s+2}^{m_j}\sum_{n=\nu_j}^{\infty}  \alpha_l^{(j)}
 e^{2\pi i \mu_j \tau}a_n(l-1-s, j)e^{2\pi i n\tau}.
 \end{eqnarray}
 
 \medskip
 Let $\tilde{\mu}_1, \hdots, \tilde{\mu}_p$ be the \emph{distinct} values among $\mu_1, \hdots, \mu_t$.
 Then we can rewrite (\ref{goodqexp1}) in the form
 \begin{eqnarray}\label{goodqexp2}
&&\ \ \ \ \ \ \ \ \ \ \ \ \  \ \ \ \ \ \   (c\tau+d)^{-k}\varphi_u(\gamma \tau)  \notag \\
&=& {\tau \choose B}\sum_{j=1}^p 
e^{2\pi i \tilde{\mu}_j \tau} g^{(j)}_B(\tau) 
+ \sum_{s=0}^{B-1} {\tau \choose s} \sum_{j=1}^p
 e^{2\pi i \tilde{\mu}_j \tau}g^{(j)}_s(\tau),
  \end{eqnarray}
where the first term in (\ref{goodqexp2}) is \emph{nonzero} and each $g^{(j)}_m(\tau)$ is a left-finite \emph{pure} $q$-series, i.e. one with only integral powers of $q$.

\bigskip
Consider the \emph{nonzero} summands
\begin{eqnarray}\label{moresum}
e^{2\pi i \tilde{\mu}_j \tau} g^{(j)}_B(\tau) &=& \sum_{n=\nu(j, B)}^{\infty} b_n(j, B)q^{n+\tilde{\mu}_j}  \\
&=&b_{\nu(j, B)}(j, B)q^{\nu(j, B)+\tilde{\mu}_j}(1+\mathrm{positive\ integral\ powers\ of}\ q) \notag
\end{eqnarray}
that occur in the first term on the right hand side of (\ref{goodqexp2}). Let $J$ be the corresponding set of indices $j$. Because the $\tilde{\mu}_j$ are distinct, there is
a \emph{unique} $j_0 \in J$ which minimizes the expression
\begin{eqnarray*}
 \nu(j, B)+\tilde{\mu}_j .
 \end{eqnarray*}
Let $J' = J\setminus{\{j_0\}}$. Hence, there is $y_0 > 0$ such that for $\Im(\tau)>y_0$ we have
\begin{eqnarray*}
\left| e^{2\pi i \tilde{\mu}_{j_0}\tau} g^{(j_0)}_B(\tau)  \right| > 2\left| \sum_{j\in J'} e^{2\pi i \tilde{\mu}_{j}\tau} g^{(j)}_B(\tau)  \right|.
\end{eqnarray*}

\medskip
\noindent
Taking into account the terms
$e^{2\pi i \tilde{\mu}_j \tau} g^{(j)}_B(\tau) $ that vanish, we obtain for $\Im(\tau)>y_0$:
\begin{eqnarray}\label{ineq2}
\left| \sum_{j=1}^p e^{2\pi i \tilde{\mu}_{j}\tau} g^{(j)}_B(\tau)  \right| &>& \left| e^{2\pi i \tilde{\mu}_{j_0}\tau} g^{(j_0)}_B(\tau)  \right| - \left| \sum_{j\in J'} e^{2\pi i \tilde{\mu}_{j}\tau} g^{(j)}_B(\tau)  \right| \notag \\
&>& 1/2 \left| e^{2\pi i \tilde{\mu}_{j_0}\tau} g^{(j_0)}_B(\tau)  \right| > 0.
\end{eqnarray}

\medskip
In (\ref{ineq2}),  for $N \in \mathbb{Z}$ we have
$\Im(\tau+N)=\Im(\tau)>y_0$. So (\ref{ineq2}) holds with $\tau$ replaced by $\tau+N$.
Because $g_B^{(j)}(\tau+N) = g_B^{(j)}(\tau)$, we see that
\begin{eqnarray}\label{ineq3}
\left| \sum_{j=1}^p e^{2\pi i \tilde{\mu}_{j}(\tau+N)} g^{(j)}_B(\tau)  \right| > 1/2 \left| e^{2\pi i \tilde{\mu}_{j_0}\tau} g^{(j_0)}_B(\tau)  \right| > 0 \ \ (N \in \mathbb{Z},\ \Im(\tau)>y_0).
\end{eqnarray}

\bigskip
At this point we return to (\ref{goodqexp2}). Replace $\tau$ by $\tau+N$ to obtain
 \begin{eqnarray}\label{goodqexp5}
&&\ \ \ \ \ \ \ \ \ \ \ \ \  \ \ \ \ \ \   (c\tau+cN+d)^{-k}\varphi_u(\gamma(\tau+N))  \notag \\
&=& {\tau+N \choose B}\sum_{j=1}^p 
e^{2\pi i \tilde{\mu}_j (\tau+N)} g^{(j)}_B(\tau) 
+ \sum_{s=0}^{B-1} {\tau+N \choose s} \sum_{j=1}^p
 e^{2\pi i \tilde{\mu}_j( \tau+N)}g^{(j)}_s(\tau).
  \end{eqnarray}

\medskip
Set
\begin{eqnarray*}
\Sigma_1(\tau)&=& \Sigma_1(\tau; \gamma) = \sum_{s=0}^{B} {\tau \choose s} \sum_{j=1}^p
 e^{2\pi i \tilde{\mu}_j( \tau)}g^{(j)}_s(\tau), \\
 \Sigma_2(\tau, N)&=&\Sigma_2(\tau, N; \gamma) = \Sigma_1(\tau+N).
\end{eqnarray*}
Thus (\ref{goodqexp5}) reads
\begin{eqnarray}\label{phisig}
\varphi_u(\gamma(\tau+N)) = (c\tau+cN+d)^k\Sigma_2(\tau, N).
\end{eqnarray}

\medskip
Next we examine the powers of $N$ that appear in $\Sigma_2(\tau, N)$. Now
\begin{eqnarray*}
{\tau +N \choose s} &=& \frac{1}{s!}(\tau+N)(\tau+N-1)\hdots(\tau+N-s+1) \\
&=& \frac{N^s}{s!} + O(N^{s-1}), \ \ N \rightarrow \infty.
\end{eqnarray*}
Therefore, the highest power of $N$ occurring with nonzero coefficient in $\Sigma_2(\tau, N)$ is $N^B$,
the coefficient in question being
\begin{eqnarray}\label{nzcoeff}
\frac{1}{B!}\sum_{j=1}^p e^{2\pi i \tilde{\mu}_j(\tau+N)}g_B^{(j)}(\tau).
\end{eqnarray}
Hence, we obtain
\begin{eqnarray}\label{nzcoeff2}
\Sigma_2(\tau, N) = \frac{N^B}{B!}\sum_{j=1}^p e^{2\pi i \tilde{\mu}_j(\tau+N)}g_B^{(j)}(\tau)
+O(N^{B-1}), \ \ N \rightarrow \infty,
\end{eqnarray}
with nonzero leading coefficient (\ref{nzcoeff}).

\bigskip
So far, the component function $\varphi_u(\tau)$ of $F(\tau)$ has been arbitrary. Now we claim that there is at least one component such that the integer $B$
occurring in (\ref{nzcoeff2}), and thereby also in (\ref{phisig}),  is equal to $M-1$. Indeed, because the first block has size $M$, the $M$th.
 component 
$\varphi_M^{(1)}(\tau)$ of $F(\tau)$ has the polynomial ${\tau \choose M-1}$ occurring in its logarithmic $q$-expansion (\ref{logqexp1}) with \emph{nontrivial} coefficient $h_0^{(1)}(\tau)$. 
Because $\rho(\gamma)$ is nonsingular, at least one row of $\rho(\gamma)$, say the $u$th., has a nonzero value $\alpha^{(1)}_{M}$ in the $M$th. column. Thanks to (\ref{conc1})
we must have $B+1=M$, as asserted.

\medskip
With $\varphi_u(\tau)$ as in the last paragraph, we have for $N \rightarrow \infty$, 
\begin{eqnarray}\label{phisig2}
&&\ \ \ \ \ \ \ \ \ \ \ \ \ \ \ \ \ \ \ \ \ \  \varphi_u(\gamma(\tau+N)) = \notag \\
&&(c\tau+cN+d)^k\left(\frac{N^{M-1}}{(M-1)!}\sum_{j=1}^p e^{2\pi i \tilde{\mu}_j(\tau+N)}g_{M-1}^{(j)}(\tau)
+O(N^{M-2})   \right).
\end{eqnarray}

\begin{lemma}\label{lemmacont} Let $\varphi_u(\tau)$ be as before, and suppose that there exists a rational number $a/c, \ ((a, c)=1, c\not= 0)$ at which $\varphi_u(\tau)$
is \emph{continuous from above}. Then $k \leq 1-M$.
\end{lemma}
\begin{pf} First note that 
\begin{eqnarray*}
\gamma(\tau+N) = \frac{a+b/(\tau+N)}{c+d(\tau+N)} \rightarrow a/c \ \mathrm{as}\ N \rightarrow \infty
\ \mathrm{within}\ \frak{H}.
\end{eqnarray*}
By the continuity assumption of the Lemma, $\varphi_u(\tau)$ remains bounded
as $N \rightarrow \infty$.
Choosing $y_0$ large enough, we see from (\ref{ineq3}) that 
\begin{eqnarray*}
\left| \sum_{j=1}^p e^{2\pi i \tilde{\mu}_{j}(\tau+N)} g^{(j)}_{M-1}(\tau)  \right| 
\end{eqnarray*}
is bounded away from $0$ as $N \rightarrow \infty$. Because $c\not= 0$, we deduce that
 the right hand side of (\ref{phisig2}) is $\approx \alpha(N)N^{k+M-1}$ with $\alpha(N)\not= 0$. If 
 $k>1-M$ this implies that the right hand side is unbounded as $N \rightarrow \infty$. This contradiction proves the Lemma. $\hfill \Box$
\end{pf}

\begin{lemma}\label{lemmageq} Let $\varphi_u(\tau)$ be as in (\ref{phisig2}), and suppose that $\varphi_u(\tau)$ is holomorphic in a region
containing $\frak{H} \cup \frak{I}$ with $\frak{I}$ a non-empty open interval in $\mathbb{R}$.
Then $k \geq 1-M$.
\end{lemma}
\begin{pf} Choose rational $a/c$ as in the last Lemma so that $a/c \in \frak{I}$. The argument of the previous Lemma
shows that the right hand side of (\ref{phisig2}) is $\approx \alpha(N)N^{k+M-1}$ with 
$\alpha(N)\not= 0$. Indeed, we easily see from (\ref{phisig2}) that $\alpha(N)$ has an upper bound independent of $N$ for $N \rightarrow \infty$. Then if $k<1-M$ the right hand side of (\ref{phisig2})
$\rightarrow 0$ as $N \rightarrow \infty$.

\medskip
On the other hand, the left hand side of (\ref{phisig2}) $\rightarrow \varphi_u(a/c)$ as 
$N \rightarrow \infty$. We conclude that $\varphi_u(a/c)=0$, and because this holds for all
rationals in $\frak{I}$ then $\varphi_u(\tau)$ is identically zero thanks to the regularity assumption on 
$\varphi_u(\tau)$. Because the components of $F(\tau)$ are linearly independent $\varphi_u(\tau)$ cannot vanish, and this contradiction shows that $k\geq 1-M$, as required. $\hfill \Box$
\end{pf}

\begin{propn}\label{propnk=1-M}Assume that the regularity assumption of Lemma \ref{lemmageq} applies to
\emph{all} components of $F(\tau)$.
Then $k=1-M$ and each component is a polynomial of degree at most $M-1$.
\end{propn}
\begin{pf} Because of the regularity assumptions of the present Proposition we may apply
Lemmas \ref{lemmacont} and \ref{lemmageq} to find that $k=1-M$. If, for some component $\varphi_u(\tau)$, the integer $B$ that occurs in (\ref{nzcoeff2}) is \emph{less} than $M-1$,
the argument of Lemma \ref{lemmacont} yields a contradiction. Since, in any case, we have
$B \leq M-1$, then in fact  $B=M-1$ for all components. As a result, (\ref{phisig2}) holds for
\emph{every} component $\varphi(\tau)$ of $F(\tau)$.

\medskip
Differentiate (\ref{goodqexp2}) (now with $B=M-1$) $M$ times and apply the well-known identity of
Bol \cite{B}: 
\begin{eqnarray*}
D^{(M)}\left((c\tau+d)^{M-1}\varphi(\gamma \tau) \right) = (c\tau+d)^{-1-M}\varphi^{(M)}(\gamma \tau).
\end{eqnarray*}
We obtain (using Leibniz's rule) for $|\tau| \rightarrow \infty$ that
\begin{eqnarray*}
&&\ \ \ \ \ \ \ \ \ \ \ \ \ \ \ \ \ \ \ \ \ \ \ \ \ \ \ \ \varphi^{(M)}(\gamma \tau) \\
&=& (c\tau+d)^{M+1}D^{(M)} \left( \sum_{s=0}^{M-1} {\tau \choose s} \sum_{j=1}^p e^{2\pi i \tilde{\mu}_j \tau}g^{(j)}_s(\tau) \right) \\
&=&(c\tau+d)^{M+1} \left(\frac{\tau^{M-1}}{(M-1)!}D^{(M)} \left(  \sum_{j=1}^p 
e^{2\pi i \tilde{\mu}_j \tau}g^{(j)}_{M-1}(\tau) \right)+O(|\tau|^{M-2}) \right).
\end{eqnarray*}
Therefore, for $N \rightarrow \infty$,
\begin{eqnarray}\label{est3}
&&\ \ \ \ \ \ \ \ \ \ \ \ \ \ \ \ \ \ \ \ \ \ \ \ \ \ \ \ \ \ \ \ \ \ \ \ \ \ \ \ \ \  \varphi^{(M)}(\gamma (\tau+N)) =  \\
&&(c\tau+cN+d)^{M+1} \left(\frac{(\tau+N)^{M-1}}{(M-1)!}D^{(M)} \left(  \sum_{j=1}^p 
e^{2\pi i \tilde{\mu}_j (\tau+N)}g^{(j)}_{M-1}(\tau) \right)+O(N^{M-2}) \right). \notag
\end{eqnarray}

\medskip
Take $c\not= 0$ with $\gamma$ as in
(\ref{abcdmatrix}) and $a/c \in \frak{I}$, and apply the regularity assumption to $\varphi(\tau)$.
Then the left hand side of (\ref{est3}) has a limit
$\varphi^{(M)}(a/c)$ for $N \rightarrow \infty$.  On the other hand, we know that  $\left|\sum_{j=1}^p 
e^{2\pi i \tilde{\mu}_j (\tau+N)}g^{(j)}_{M-1}(\tau)\right|$ is bounded away from zero. So if
$D^{(M)} \left(  \sum_{j=1}^p 
e^{2\pi i \tilde{\mu}_j (\tau+N)}g^{(j)}_{M-1}(\tau) \right)$ does not vanish identically then
the right hand side of (\ref{est3}) is $\approx \alpha(N)N^{2M}$ for
$N \rightarrow \infty$, with $\alpha(N)$ bounded away from zero. This contradiction shows that in fact
\begin{eqnarray*}
D^{(M)} \left(  \sum_{j=1}^p 
e^{2\pi i \tilde{\mu}_j (\tau+N)}g^{(j)}_{M-1}(\tau) \right) \equiv 0
\end{eqnarray*}

\medskip
From (\ref{moresum}), we have
\begin{eqnarray*}
&&\ \ \ \ \ \ \ \ \  D^{(M)}\left( \sum_{j=1}^p e^{2\pi i \tilde{\mu}_j (\tau+N)} g^{(j)}_{M-1}(\tau) \right) \\
&=&\sum_{j=1}^pe^{2\pi i \tilde{\mu}_jN}
 \sum_{n=\nu(j, M-1)}^{\infty} b_n(j, M-1)(2\pi i (n+\tilde{\mu}_j))^Mq^{n+\tilde{\mu}_j},
\end{eqnarray*}
so that
\begin{eqnarray*}
\sum_{j=1}^p e^{2\pi i \tilde{\mu}_jN} b_n(j, M-1)(n+\tilde{\mu}_j)^M = 0 \ \  \ \ (n \geq \nu(j, M-1)).
\end{eqnarray*}
This implies that $b_n(j, M-1)=0$ whenever $n+\tilde{\mu}_j\not= 0$.
Because $b_{\nu(j, M-1)}(j, M-1)\not= 0$,  we must have 
\begin{eqnarray}\label{muvan}
\tilde{\mu}_j = 0 \ \ \ (1 \leq j \leq p), 
\end{eqnarray}
which amounts to the assertion that all $\mu_j=0\ (1 \leq j \leq t)$. Moreover
$b_n(j, M-1)=0$ for $n\not= 0$, so that
\begin{eqnarray*}
g_{M-1}^{(j)}(\tau) = b_0(j, M-1)
\end{eqnarray*}
is  constant. Now (\ref{goodqexp2}) reads
 \begin{eqnarray}\label{goodqexp21/2}
  (c\tau+d)^{M-1}\varphi(\gamma \tau) 
= {\tau \choose M-1}\sum_{j=1}^p b_0(j, M-1)
+ \sum_{s=0}^{M-2} {\tau \choose s} \sum_{j=1}^p g^{(j)}_s(\tau).
  \end{eqnarray}

\medskip
We now repeat the argument $M-1$ times, starting with (\ref{goodqexp21/2}) in place of (\ref{goodqexp2}). We end up with an identity of the form
 \begin{eqnarray*}
  (c\tau+d)^{M-1}\varphi(\gamma \tau) 
=  \sum_{s=0}^{M-1} {\tau \choose s} \sum_{j=1}^p b_0(j, s),
  \end{eqnarray*}
where of course the right hand side is a polynomial $p(\tau)$ of degree at most $M-1$.
Then
\begin{eqnarray*}
\varphi(\tau) &=& (c\gamma^{-1}\tau+d)^{1-M}p\left( \frac{\ \ d\tau-b}{-c\tau+a}\right) \\
&=& (c\tau+d)^{M-1}p\left( \frac{\ \ d\tau-b}{-c\tau+a}\right) 
\end{eqnarray*}
is itself a polynomial of degree at most $M-1$. This completes the proof of  Proposition \ref{propnk=1-M}.
$\hfill \Box$
\end{pf}

\bigskip
It is now easy to complete the proof of Theorem 2. It is only necessary to establish the implication
(a) $\Rightarrow$ (b). Assuming (a)  means that  Proposition \ref{propnk=1-M} is applicable, so that we have $k=1-M$ and the components of $F(\tau)$ are polynomials of degree at most $M-1$. Because the components are linearly independent, it must be the case that the maximal block size $M$ is equal to 
 the dimension $p$ of the representation $\rho$. Thus $k=1-p$, and the component functions  \emph{span} the space of polynomials of degree at most $p-1$. The fact that $(\rho, F)$ is equivalent to
 $(\sigma, C)$ then follows from the discussion in Section \ref{S2}.


\begin{thebibliography}{BPZ}
   
   \bibitem[B]{B} Bol, G., Invarianten linearer Differentialgleichungen, Abh. Math. Sem. Univ. Hamburg \textbf{16} Nos. 3-4 (1949), 1-28.
 
 \bibitem[H]{H} Hille, E., \textit{Ordinary Differential Equations in the
 Complex Domain}, Dover Publications, New York, 1976. 
 
 \bibitem[KM1]{KM1} Knopp, M., and Mason, G.,  On vector-valued modular forms and their Fourier coefficients,
 Acta Arith. \textbf{110.2} (2003), 117-124.
 
 \bibitem[KM2]{KM2} Knopp, M. and Mason, G., Vector-valued modular forms and Poincar\'{e} series, 
 Ill. J. Math. \textbf{48} No. 4 (2004),  1345-1366.
 
 
  
 \bibitem[KM3]{KM3} Knopp, M. and Mason, G., Logarithmic vector-valued modular forms, to appear in
 Acta Arithmetica.
 
 
  \bibitem[M1]{M1} Mason, G., Vector-valued modular forms and linear differential equations, 
Int. J. Numb. Th. \textbf{3} No. 3 (2007), 1-14.


 
  \end{thebibliography}
\end{document}